\documentclass[11pt,reqno]{amsart}
\usepackage[english]{babel}
\usepackage[latin1]{inputenc}
\usepackage{fancyhdr}
\usepackage{indentfirst}
\usepackage{graphicx}
\usepackage{float}
\usepackage{caption}
\usepackage{newlfont} 
\usepackage{amssymb}
\usepackage{amsmath}
\usepackage{mathtools}
\usepackage{latexsym}
\usepackage{amsthm}
\usepackage{enumitem}
\usepackage{esint}
\usepackage{amsfonts}
\usepackage{amscd,amsxtra}
\usepackage{color}

\usepackage[english]{babel}
\usepackage[latin1]{inputenc}
\usepackage{fancyhdr}
\usepackage{indentfirst}
\usepackage{graphicx}
\usepackage{float}
\usepackage{caption}
\usepackage{newlfont} 
\usepackage{amssymb}
\usepackage{amsmath}
\usepackage{mathtools}
\usepackage{latexsym}
\usepackage{amsthm}
\usepackage{enumitem}
\usepackage{amsfonts}
\usepackage{amscd,amsxtra}

\newcommand{\scal}[2]{\langle {#1} , {#2}\rangle}

\newcommand{\divg}{\mathrm{div}_\mathbb{G}\,}

\newtheorem{thm}{Theorem}[section]
\newtheorem{cor}[thm]{Corollary}

\newtheorem{prop}[thm]{Proposition}
\theoremstyle{definition}
\newtheorem{defn}[thm]{Definition}
\theoremstyle{remark}
\newtheorem{rem}[thm]{Remark}
\numberwithin{equation}{section}

\makeatletter
\def\cleardoublepage{\clearpage\if@twoside \ifodd\c@page\else
	\hbox{}
	\thispagestyle{empty}
	\newpage
	\if@twocolumn\hbox{}\newpage\fi\fi\fi}
\makeatother
\title{Alt-Caffarelli-Friedman monotonicity formula and mean value properties in Carnot groups with applications}
\author{Fausto Ferrari}
\address{Fausto Ferrari: Dipartimento di Matematica\\ Universit\`a di Bologna\\ Piazza di Porta S.Donato 5\\ 40126, Bologna-Italy}
\email{fausto.ferrari@unibo.it }
\author{Nicol\`o Forcillo}
\address{Nicol\`o Forcillo: Dipartimento di Matematica\\ Universit\`a di Roma Tor\linebreak Vergata\\ Via della Ricerca Scientifica 1 \\ 00133, Roma-Italy}
\email{nikforc6@gmail.com}
\thanks{F.F.  and N.F. are partially supported by INDAM-GNAMPA-2019 project: {\it Propriet\`a di regolarit\`a delle soluzioni viscose con applicazioni a problemi di frontiera libera.}}
\date{\today}
\begin{document}
	\begin{abstract}
		In this paper we provide a different approach to the Alt-Caffarelli-Friedman monotonicity formula, reducing the problem to test the monotone increasing behavior of the mean value of a function involving the norm of the gradient. 
		  In particular,  we show that our argument holds in the general framework of Carnot groups.
	\end{abstract}
	\maketitle

    \tableofcontents	
	\section{Introduction}
	In this paper, we discuss sufficient conditions under which an Alt-Caffarelli-Friedman-type monotonicity formula holds for functions satisfying homogeneous second order elliptic equations in Carnot groups. 
	
	In addition, we provide a different approach to that type of result that changes the classical strategy of facing the problem and becomes useful for obtaining counterexamples in noncommutative structures. We refer to \cite{ACF} and \cite{CS} for the original proofs in the Euclidean setting, \cite{FeFo} for an overview, and \cite{FeFo2} for a counterexample in the Heisenberg group.
	
	We have been motivated in this research by the project of extending some regularity results for free boundary problems already known for linear and nonlinear elliptic operators, see \cite{FeLeSa}.  
	
	Our strategy is mainly based on the geometric structure underlying the
	functions we are dealing with. The main tool is the existence of a
	group of dilations under which the solutions remain in the same class. We can refer to this as an invariant by scaling property. Moreover, the existence of a fundamental solution associated with the operator and, as a byproduct, the validity of intrinsic mean value properties for intrinsic harmonic functions are crucial.
	
    The proof we are going to recall in the next section could be considered a folk approach to the problem of finding mean formulas. Indeed, its trace can be found in the potential theory literature for a long time. It was first applied to the Laplace operator and then adapted to different frameworks, for instance, taking into account the geometry associated with the operator.

		An essential remark, in the problem of finding mean formulas, was made when it was noticed that the set on which it is useful to apply the mean formula has to be determined by the level surfaces of the fundamental solution of the operator under consideration. As far as we know, this idea might be consider a folk knowledge. The papers in which it was first explicitly used should be \cite{Weber},  \cite{Pini_1}, \cite{Pini_2}, \cite{Fulks}, \cite{Kupcov}, \cite{Kupcov2}, \cite{Kupcov3}, \cite{Kuptsov} and, more recently, it was applied in \cite{Fabes_Garofalo} in the parabolic setting. In  \cite{Citti_Garofalo_Lanconelli},  the approach was applied to degenerate elliptic operators, determining a new natural geometry. We quote here \cite{Hormander}, \cite{Folland_Stein}, and \cite{BLU} for a comprehensive discussion in Carnot groups. Also, in \cite{Vesely} an interesting  historical  review of the mean value properties can be found. In this paper, we mainly use the results in \cite{Citti_Garofalo_Lanconelli}, see also \cite{BLU}, where the theory about mean value properties was completely developed and fixed.

		 In order to introduce our main results, we would need to discuss many definitions that we have postponed to Section \ref{dueF} for organizing reasons. Let us simply recall that $\mathcal{L}$ denotes a second order operator satisfying some suitable properties, such as the class of sub-Laplacians on Carnot groups. People unfamiliar with this stuff may think that  $\mathcal{L}$ is the usual Laplace operator on the Euclidean space $\mathbb{R}^N$ and $\Omega\subset \mathbb{R}^N$ is an open connected set. Moreover, given a function $u$ in a suitable Sobolev space, we write the following functional:
		\begin{equation}\label{moneqH01}
		J_{u,x_0}^{\mathbb{G}}(r)\coloneqq \frac{1}{r^4}\int_{D(x_0,r)}|\nabla_{\mathbb{G}} u^+|^2\Gamma(x_0,\xi)\hspace{0.075cm}d\xi\int_{D(x_0,r)}|\nabla_{\mathbb{G}} u^-|^2\Gamma(x_0,\xi)\hspace{0.075cm}d\xi,
	\end{equation}
		where $\Gamma$ denotes the fundamental solution associated with the operator $\mathcal{L},$ $D(x_0,r)$ is the appropriate ball centered at $x_0$ of radius $r,$ and $\mathbb{G}$ is a Carnot group, see Section \ref{dueF} for the precise definitions. We recall that \linebreak $u^+\coloneqq\sup\{u,0\}$, $u^-\coloneqq\sup\{-u,0\}$. For instance, if $\mathbb{G}=\mathbb{R}^N$, $N\in \mathbb{N}$, endowed with the usual commutative inner law, $\mathcal{L}=\Delta$, $\nabla_{\mathbb{G}}=\nabla$,  \linebreak $\Gamma (x,y)=c_N|x-y|^{2-N}$, $D(x_0,r)=B(x_0,r)=\{x\in \mathbb{R}^N:\hspace{0.2cm} |x-x_0|_{\mathbb{R}^N}<r\},$ and $|\cdot|_{\mathbb{R}^N}$ denotes the usual Euclidean norm in $\mathbb{R}^N$, then
		\begin{equation}\label{moneqH0}
		J_{u, x_0}^{\mathbb{R}^N}(r)\coloneqq c_N^2\frac{1}{r^4}\int_{B(x_0,r)}\frac{|\nabla u^+|^2}{|\xi-x_0|_{\mathbb{R}^N}^{N-2}}\hspace{0.075cm}d\xi\int_{B(x_0,r)}\frac{|\nabla u^-|^2}{|\xi-x_0|_{\mathbb{R}^N}^{N-2}}\hspace{0.075cm}d\xi
	\end{equation}
		becomes, up to a multiplicative constant, the classical functional associated with the Alt-Caffarelli-Friedman formula in the Euclidean framework, see \cite{ACF}. In \cite{ACF}, Alt, Caffarelli and Friedman prove that if $\Delta u^{\pm}\geq 0$ in $\Omega$ and $x_0\in\partial \Omega^+(u)$, then $J_{u, x_0}^{\mathbb{R}^N}$ is monotone increasing in a right neighborhood of $0$.
		
		We state here a stronger version of our main result, and we refer to Section \ref{treF} for a more general Theorem \ref{maintheoremgeneral}, including some technical definitions on the mean formulas in Carnot groups, see \cite{BLU}, which we prefer to avoid introducing here for simplifying the exposition.

	Mainly following the notation of Section 5.5 in \cite{BLU}, for every $x_0\in \mathbb{G}$ and $r>0$, we define the intrinsic ball of center $x_0$ and radius $r$ as  $D(x_0,r):=\{y\in \mathbb{G}:\hspace{0.2cm}|x_0^{-1}\circ y|<r\},$ where $|\cdot |$ denotes a suitable norm on the group $\mathbb{G}$. Moreover, for every continuous function $u\in C(\mathbb{G})$ and $r>0,$ we denote the mean value associated with $u$ and $r$ at $x_0$ in the Carnot group $\mathbb{G}$ by:
	\begin{equation}\label{averageformulaint}
	\mathcal{M}_r(u)(x_0):=\int_{\partial D(x_0,r)}u(y)\mathcal{K}(x_0^{-1}\circ y)d\sigma(y),
	\end{equation}
	where $\mathcal{K}$ is a positive kernel that we recall in Section \ref{dueF}, together with further details about previous definitions. We set $\Omega^+(u):=\{x\in \Omega:\:u(x)>0\}$ and $\Omega^-(u):=\{x\in \Omega:\:u(x)\leq 0\}^o$.

 		\begin{thm}\label{maintheorem}
		Let $\Omega\subset \mathbb{G}$ be an open set in a Carnot group $\mathbb{G}.$ Let $u\in C(\Omega)$ such that $|\nabla_{\mathbb{G}} u|^2\in L^1_{\mathrm{loc}}(\Omega)$.
		For every $x_0\in \partial \{u>0\}\cap \Omega$, if 
		\begin{equation}\label{formfond}r\to\mathcal{M}_r\left(\left(\frac{|\nabla_{\mathbb{G}} u^{\pm}(x_0\circ\cdot)|}{|\nabla_{\mathbb{G}}|\cdot||}\right)^2\right)(0)\end{equation}  are monotone increasing almost everywhere in a right neighborhood of $0$, then the function $J_{u, x_0}^{\mathbb{G}}$ is monotone increasing in a right neighborhood of $0$.
		\end{thm}
		Therefore, when $\mathbb{G}={\mathbb{R}^N}$, condition \eqref{formfond} reduces to verify the local increasing behavior almost everywhere of
		$$
		r\to\mathcal{M}_r\left(|\nabla u^{\pm}(x_0+\cdot)|^2\right)(0),
		$$
		because $|\nabla |x||=1$. In particular, as a consequence of our Theorem \ref{maintheorem}, the following result holds. 
		\begin{cor}\label{yes}
		Let $u\in C^2(\Omega)$ be harmonic in $\Omega$. For every $x_0\in \Omega$, 
		\begin{equation}\label{moneqH0E}
		r\to \frac{c_N}{r^2}\int_{B(0,r)}\frac{|\nabla u(x_0+\xi)|^2}{|\xi|_{\mathbb{R}^N}^{N-2}}\hspace{0.075cm}d\xi	\end{equation}
		is monotone increasing in a right neighborhood of $0$.
		\end{cor}

		
		\begin{rem}
		Our results appear useful in some cases. Nevertheless, it cannot be considered an alternative proof of the celebrated Alt-Caffarelli-Friedman monotonicity theorem. In fact, 
		our hypotheses apply to the case in which
	  $\Delta u=0$ in $\Omega^+(u)$ and $\Delta u=0$ in $\Omega^-(u)$, that is in the same assumptions of the Alt-Caffarelli-Friedman monotonicity formula. However, despite
		$u^+$ and $u^-$ are subharmonic, we cannot deduce that $\xi\to|\nabla u^{\pm}(x_0+\xi)|^2$ are subharmonic. Indeed, in general, $\xi\to|\nabla u^{\pm}(x_0+\xi)|^2$ are not upper semicontinuous, so they cannot be subharmonic, and we cannot deduce that $r\to\mathcal{M}_r\left(|\nabla u^{\pm}(x_0+\cdot)|^2\right)(0)$ are monotone increasing. In any case, in some particular situations, by knowing the relationship between $
		r\to\mathcal{M}_r\left(|\nabla u^{+}(x_0+\cdot)|^2\right)(0) 
		$
		and 
		$
		r\to\mathcal{M}_r\left(|\nabla u^{-}(x_0+\cdot)|^2\right)(0),
		$
		we may deduce some results about the monotonicity of $J_{u, x_0}^{\mathbb{R}^N}$.
		\end{rem}
		The last argument may be used in Carnot groups for constructing counterexamples to the increasing monotone behavior of $J_{u,x_0}^{\mathbb{G}}$. Indeed, in this paper, in addition to the characterization via mean formulas, we obtain new counterexamples to the monotonicity of $J_{u, x_0}^{\mathbb{H}^1}$ in the Heisenberg group $\mathbb{H}^1,$ see also  \cite{FeFo2}. Specifically, we obtain the following result.		

		\begin{cor}\label{maintheorem2}
		In the Heisenberg group $\mathbb{H}^1,$ there exist harmonic functions and points  $x_0\in \partial \{u>0\}\cap \Omega$ such that 
		  $J_{u, x_0}^{\mathbb{H}^1}$ is monotone decreasing in a right neighborhood of $0$. In particular,
		  for every $c_1,c_2\in \mathbb{R}$ such that $c_1^2+c_2^2>0,$ the function $$u(x,y,t)=c_1x+c_2y+3t(c_2x-c_1y)-2(c_1x^3+c_2y^3)$$ is harmonic in $\mathbb{H}^1,$ $0\in \partial \{u>0\},$ and 
		  $$r\to J_{u, 0}^{\mathbb{G}}(r)$$
		  is monotone decreasing in a right neighborhood of $0$.
		\end{cor}
		
		The paper is organized as follows. In Section \ref{unoF}, we start by recalling the classical approach described in \cite{GT} to determine mean formulas. In Section \ref{dueF}, we introduce the basic notation of Carnot groups. In Section \ref{treF}, we show the proofs of our main results and corollaries.

\section{The well-known key idea for obtaining a mean value formula}\label{unoF}	
	In this section, we recall how to obtain a mean value formula in a standard way. Let $\mathcal{O}\subset \mathbb{R}^n$ be an open bounded connected set. For every open set $\Omega\subset \mathcal{O},$ endowed with a smooth boundary, and for every $y\in \Omega,$ we consider a function $u$ satisfying
	$$
	\mathcal{L}u:=\mbox{div}(A(x)\nabla u(x))=f(x)\quad\mbox{in }\mathcal{O},
	$$
	where $A\geq 0$ is a matrix with smooth coefficients. If the fundamental solution $\Gamma_A$ exists for $\mathcal{L}$, then, denoting by $D(y)\subset\Omega$ a smooth neighborhood of $y$, it holds
	\begin{equation*}\label{unoFa}\begin{split}
	&\int_{\Omega\setminus D(y)} \Gamma_A(y,x)f(x)\hspace{0.075cm}dx=\int_{\Omega\setminus D(y)} \left(\Gamma_A(y,x)\mathcal{L}u(x)-u(x)\mathcal{L}\Gamma_A(y,x)\right) dx\\
	&=\int_{\partial (\Omega\setminus D(y))}\left(\langle A(x)\nabla u(x),n\rangle\hspace{0.025cm}\Gamma_A(y,x)-\langle A(x)\nabla \Gamma_A(y,x),n\rangle\hspace{0.025cm} u(x)\right)d\sigma(x).
	\end{split}
	\end{equation*}
	In particular, if $\Omega =\Omega_R(y):=\left\{x\in \mathcal{O}:\:\:\Gamma_{A}(x,y)>\frac{1}{R}\right\},$ $R>0,$ and we take $D(y)=\Omega_r(y),$ $0<r<R,$ then
	\begin{equation}\label{treFa0}\begin{split}
	&\int_{\Omega_R(y)\setminus \Omega_r(y)} \Gamma_A(y,x)f(x)\hspace{0.05cm}dx=\\
	&=\int_{\partial (\Omega_R(y)\setminus \Omega_r(y))}\langle A(x)\nabla u(x),n\rangle\hspace{0.025cm}\Gamma_A(y,x)\hspace{0.075cm}d\sigma(x)\\
	&-\int_{\partial (\Omega_R(y)\setminus \Omega_r(y))}\langle A(x)\nabla \Gamma_A(y,x),n\rangle\hspace{0.025cm}u(x)\hspace{0.075cm}d\sigma(x)\\
	&=\frac{1}{R}\int_{\partial\Omega_R(y)}\langle A(x)\nabla u(x),n\rangle\hspace{0.05cm} d\sigma(x)-\frac{1}{r}\int_{\partial \Omega_r(y)}\langle A(x)\nabla u(x),n\rangle\hspace{0.05cm} d\sigma(x)\\
	&-\int_{\partial\Omega_R(y)}\langle A(x)\nabla \Gamma_A(y,x),n\rangle\hspace{0.025cm} u(x)\hspace{0.075cm}d\sigma(x)\\
	&+\int_{\partial \Omega_r(y)}\langle A(x)\nabla \Gamma_A(y,x),n\rangle\hspace{0.025cm} u(x)\hspace{0.075cm}d\sigma(x).\\
	\end{split}
	\end{equation}

On the other hand, we know that $|\langle A(x)\nabla u(x),n\rangle|\leq \|A\|_{L^\infty(\mathcal{O})}\|\nabla u\|_{L^\infty(\mathcal{O})}.$ Thus, if $\frac{|\partial \Omega_r(y)|_{n-1}}{r}\to 0$ as $r\to 0,$ we get
$$
\frac{1}{r}\int_{\partial \Omega_r(y)}\langle A(x)\nabla u(x),n\rangle \hspace{0.075cm}d\sigma(x)\to 0\quad\text{as }r\to 0.
$$
Moreover, supposing there exists a change of variables $T_{y,r}$ such that
\begin{equation*}\begin{split}
&\int_{\partial \Omega_r(y)}\langle A(x)\nabla \Gamma_A(y,x),n\rangle\hspace{0.025cm} u(x)\hspace{0.075cm}d\sigma(x)\\
&=r^{Q-1}\int_{\partial \Omega_1(0)}\langle A(T_{y,r}(\xi))(\nabla \Gamma_A)(y, T_{y,r}(\xi)),n\rangle\hspace{0.025cm} u(T_{y,r}(\xi))\hspace{0.075cm}d\sigma(\xi)
\end{split}
\end{equation*}
and $$r^{Q-1}\int_{\partial \Omega_1(0)}\langle A(T_{y,r}(\xi))(\nabla \Gamma_A)(y,T_{y,r}(\xi)),n\rangle\hspace{0.025cm} u(T_{y,r}(\xi))\hspace{0.075cm}d\sigma(\xi)\to c_n u(y)$$ as $r\to 0$, it holds, according to \eqref{treFa0},
\begin{equation}\label{treFa2}\begin{split}
	&\int_{\Omega_R(y)\setminus \Omega_r(y)} \Gamma_A(y,x)f(x)dx=-\int_{\partial\Omega_R(y)}\langle A(x)\nabla \Gamma_A(y,x),n\rangle\hspace{0.025cm} u(x)\hspace{0.075cm}d\sigma(x)\\
	&+c_nu(y)
	+\frac{1}{R}\int_{\partial\Omega_R(y)}\langle A(x)\nabla u(x),n\rangle \hspace{0.075cm}d\sigma(x).
	\end{split}
	\end{equation}
	Assuming now that $h$ satisfies
	\begin{equation}\label{def-h}
	\begin{cases}
	\mathcal{L}h=0&\text{in }\Omega_R(y),\\
	h=\Gamma_A(\cdot,y)&\text{on }\partial \Omega_R(y),
	\end{cases}
	\end{equation}
     we achieve
\begin{equation}\label{unobisFa}\begin{split}
	&\int_{\Omega_R(y)} h(x)\mathcal{L}u(x)\hspace{0.075cm}dx=\int_{\Omega_R(y)} \left(h(x)\mathcal{L}u(x) -u(x)\mathcal{L}h(x)\right)dx\\
	&=\int_{\partial (\Omega_R(y))}\left(\langle A(x)\nabla u(x),n\rangle\hspace{0.025cm} h(x)-\langle A(x)\nabla h(x),n\rangle\hspace{0.025cm} u(x)\right)d\sigma(x).
	\end{split}
	\end{equation}
Hence, subtracting \eqref{treFa2} to \eqref{unobisFa} term by term, we find
\begin{equation*}\label{quattroFa1}\begin{split}
	&\int_{\Omega_R(y)} (h(x)-\Gamma_A(y,x))f(x)\hspace{0.075cm}dx
	\\&=\int_{\partial\Omega_R(y)}\langle A(x)\nabla (-h(x)+\Gamma_A(x,y)),n\rangle\hspace{0.025cm} u(x)\hspace{0.075cm}d\sigma(x)\\
	&-c_nu(y)
	+\int_{\partial\Omega_R(y)}\langle A(x)\nabla u(x),n\rangle\hspace{0.025cm} (-\Gamma_A(x,y)+h) \hspace{0.075cm}d\sigma(x).
	\end{split}
	\end{equation*}
Next, recalling that $\Gamma_A(x,y)-h=0$ on $\partial \Omega_R(y)$ in view of \eqref{def-h}, we obtain the following representation formula
\begin{equation*}\label{quattroFa2}\begin{split}
	c_nu(y)&=\int_{\partial\Omega_R(y)}\langle A(x)\nabla (\Gamma_A(x,y)-h(x)),n\rangle\hspace{0.025cm} u(x)\hspace{0.075cm}d\sigma(x)\\
	&-\int_{\Omega_R(y)} (h(x)-\Gamma_A(x,y))f(x)\hspace{0.075cm}dx,
	\end{split}
	\end{equation*}
	that is, assuming $c_n=1$ possibly adjusting the constant in $\Gamma_A,$
	\begin{equation*}\label{quattroFa222}\begin{split}
	u(y)&=\int_{\partial\Omega_R(y)}\langle A(x)\nabla G_{\Omega_R(y),A(x)}(x),n\rangle\hspace{0.025cm} u(x)\hspace{0.075cm}d\sigma(x)\\
	&+\int_{\Omega_R(y)}G_{\Omega_R(y),A(x)}(x)f(x)\hspace{0.075cm}dx,
	\end{split}
	\end{equation*}
	where $G_{\Omega_R(y),A(x)}$ is the Green function associated with the set $\Omega_R(y)$ evaluated at the pole $y$.
	This argument, when $f=0$, yields the following mean formula:
	\begin{equation*}\label{cinqueFa}\begin{split}
	u(y)=\int_{\partial\Omega_R(y)}\langle A(x)\nabla G_{\Omega_R(y),A(x)}(x),n\rangle\hspace{0.025cm} u(x)\hspace{0.075cm}d\sigma(x).
	\end{split}
	\end{equation*}

	Exploiting the previous strategy, with the suitable matrix $A,$ it is possible to tailor mean value formulas in Carnot groups associated with the sub-Laplacian $\Delta_{\mathbb{G}}$. A detailed overview of the tools to study mean value formulas in Carnot groups can be found in \cite{BLU}. Nevertheless, in the next section, we recall the main definitions and results we use in our proofs. 
	
	\section{Main notation of Carnot groups and mean value formulas}\label{dueF}
	
	A simply connected stratified nilpotent 
 Lie group $(\mathbb{G},\circ)$  (in general noncom\-mu\-ta\-tive) is said a {\it Carnot group  of
step $\kappa$} if its Lie algebra
${\mathfrak{g}}$  admits a {\it step $\kappa$ stratification}, i.e.,
there exist linear subspaces $V_1,...,V_\kappa$ such that
\begin{equation}\label{stratificazione}
{\mathfrak{g}}=V_1\oplus...\oplus V_\kappa,\quad [V_1,V_i]=V_{i+1},\quad
V_\kappa\neq\{0\},\quad V_i=\{0\}{\textrm{ if }} i>\kappa,
\end{equation}
where $[V_1,V_i]$ is the subspace of ${\mathfrak{g}}$ generated by
the commutators $[X,Y]$ with $X\in V_1$ and $Y\in V_i$. 
The first layer $V_1$, the so-called 
 {\sl horizontal layer},  plays a
key role in the theory, since it generates  $\mathfrak g$ by commutation. Besides \textit{\textit{\textit{\textit{}}}}Euclidean spaces, the only abelian Carnot groups, typical examples are Heisenberg groups and upper triangular groups. For a general introduction to Carnot groups, from our point of view, and further examples, we refer, e.g., to \cite{BLU}, \cite{Folland_Stein}, and \cite{stein}.

Let us consider a Carnot group $\mathbb{G}.$ Set
$m_i\coloneqq\dim(V_i)$, for $i=1,\dots,\kappa,$ and $h_i=m_1+\dots +m_i$, so that $h_\kappa=N$. For the sake of simplicity, we write 
$m\coloneqq m_1$.
We also denote by $Q$ the {\sl homogeneous dimension} of $\mathbb{G}$, i.e., we define
 $$
 Q:=\sum_{i=1}^{\kappa} i \dim(V_i).
 $$
We recall that, if $e$ is the unit element of $(\mathbb{G},\circ)$,  the map $X\to X(e)$,
which associates a left-invariant vector field $X$ with its value at $e$, is an
isomorphism from $\mathfrak g$ to $T\mathbb{G}_e$, identified with $\mathbb{R}^N$.
From now on, we shall systematically use these identifications. 

We now choose  a basis $e_1,\dots,e_N$ of
$\mathbb{R}^N$ adapted to the stratification of $\mathfrak g$, namely 
$$e_{h_{j-1}+1},\dots,e_{h_j}\;\text {is a basis of}\; V_j,\;\text{
for each}\; j=1,\dots, \kappa.$$
Moreover, let $X=\{X_1,\dots,X_{N}\}$ be the family
of left-invariant vector fields such 
that 
$X_i(e)=e_i$, $i=1,\dots,N$. 
The sub-bundle of the tangent bundle $T\mathbb{G}$ spanned by the
vector fields $X_1,\dots,X_{h_1}$ plays a crucial
role in the theory, and it is called the {\it horizontal bundle}
$H\mathbb{G}.$ The fibers of $H\mathbb{G}$ are $$ H\mathbb{G}_x=\mbox{span
}\{X_1(x),\dots,X_{h_1}(x)\},\qquad x\in\mathbb{G} .$$ We can
endow each fiber $H\mathbb{G}_x$ with a
scalar product $\scal{\cdot}{\cdot}_{x}$ and a norm
$\vert\cdot\vert_x$ that make the basis $X_1(x),\dots,X_{h_1}(x)$
an orthonormal basis. Precisely, if
$v=\sum_{i=1}^{h_1}v_iX_i(x)=(v_1,\dots,v_{h_1})$ and
$w=\sum_{i=1}^{h_1}w_iX_i(x)=(w_1,\dots,w_{h_1})$ are in $H\mathbb{G}_x$, then
$\scal{v}{w}_{x}:=\sum_{j=1}^{h_1} v_jw_j$ and $
|v|_x^2:=\scal{v}{v}_{x}$. The sections of $H\mathbb{G}$ are called {\it horizontal sections} and are the so-called {\sl horizontal vector fields}. A vector of $H\mathbb{G}_x$ is a {\it horizontal vector}, while any vector in $T\mathbb{G}_x$ not horizontal is vertical.
Each horizontal section is identified by its canonical coordinates
with respect to the moving frame $X_1(x),\dots,X_{h_1}(x)$. In this
way, a horizontal section $\phi$ is described by a function
$\phi =(\phi_1,\dots,\phi_{h_1})
:\mathbb{R}^{N} \rightarrow\mathbb{R}^{h_1}$. Moreover, a Carnot group $\mathbb{G}$ can always be identified, through exponential coordinates,
with the Euclidean space $(\mathbb{R}^N, \cdot)$,
where $N$ is the dimension of ${\mathfrak{g}}$, endowed with a suitable
group operation. The explicit
expression of the group operation $\cdot$ is determined by the
Campbell-Hausdorff formula.

For
any $x\in\mathbb{G}$, the {\it (left) translation} $\tau_x:\mathbb{G}\to\mathbb{G}$ is defined
as $$ z\mapsto\tau_x z\coloneqq x\circ z. $$ For any $\lambda >0$, the
{\it dilation} $\delta_\lambda:\mathbb{G}\to\mathbb{G}$ has the expression
\begin{equation}\label{dilatazioni}
\delta_\lambda(x_1,...,x_N)=
(\lambda^{d_1}x_1,...,\lambda^{d_N}x_N),
\end{equation} where $d_i\in\mathbb{N}$ is called {\it homogeneity of
the variable} $x_i$ in
$\mathbb{G}$ (see \cite{Folland_Stein}, Chapter 1) and is defined as
\begin{equation}\label{omogeneita2}
d_j=i \quad\text {whenever}\; h_{i-1}+1\leq j\leq h_{i}.
\end{equation}
Hence, it holds $1=d_1=...=d_{h_1}<
d_{{h_1}+1}=2\leq...\leq d_N=\kappa.$ Throughout this paper, by $\mathbb{G}$-homogeneity we mean homogeneity with respect
to group dilations $\delta_\lambda$ (see again Chapter 1 in \cite{Folland_Stein}). The Haar measure of $\mathbb{G}=(\mathbb{R}^N,\circ)$ is the Lebesgue measure
in $\mathbb{R}^N$. 
If $A\subset \mathbb{G}$ is $ L$-measurable, we
write $|A|$ to denote its Lebesgue measure. We notice that the
Lebesgue measure is invariant under the map $y\to y^{-1}$.
In addition, if $m\ge 0$, we denote the $m$-dimensional Hausdorff measure, obtained from the
Euclidean distance in $\mathbb{R}^N\simeq \mathbb{G},$ by $\mathcal
H^m.$


The following result is contained in  \cite{Folland_Stein}, see Proposition 1.26.
\begin{prop}\label{campi omogenei0}
If $j=1,\dots,m$, the vector fields $X_j$ have polynomial coefficients and the
form
\begin{equation*}\label{campi omogenei}
X_j(x)=\partial_j+\sum_{k\geq 1,\hspace{0.05cm} d_k>1} p_{j,k}(x)\partial_k, \quad
\end{equation*} where
$p_{j,k}$ are $\mathbb{G}$-homogeneous polynomials of degree $d_k-1.$
\end{prop}

We now fix a basis $X_1,\dots,X_{m}$ of the horizontal layer. For any function $f:\mathbb{G}\to \mathbb{R}$ for which the partial derivatives
$X_jf,$ $j=1,\ldots, m,$ exist, we
define the horizontal gradient of $f$, denoted by
$\nabla_{\mathbb{G}}f$, as the horizontal section
\begin{equation*}
\nabla_{\mathbb{G}}f:=\sum_{i=1}^{m}(X_if)X_i,
\end{equation*}
whose coordinates are $(X_1f,...,X_{m}f)$. Moreover, if
$\phi=(\phi_1,\dots,\phi_{m})$ is a horizontal section such
that $X_j\phi_j\in L^1_{\mathrm{loc}}(\mathbb{G}),$ $j=1,\dots,m$, we define $\divg
(\phi)$ as the real-valued function
\begin{equation*}
\divg(\phi):=-\sum_{j=1}^{m}X_j^*\phi_j=\sum_{j=1}^{m}X_j\phi_j.
\end{equation*}
According to \cite{Folland_Stein}, we also adopt the following multi-index notation for higher-order derivatives. If $I =
(i_1,\dots,i_{n})$ is a multi--index, we set  
\[X^I\coloneqq X_1^{i_1}\cdots
X_{n}^{i_{n}}.\]
By the Poincar\'e--Birkhoff--Witt theorem
(see, e.g., \cite{bourbaki}, I.2.7), the differential operators $X^I$ form a basis for the algebra of left-invariant
differential operators in $\mathbb{G}$. Furthermore, we define 
$|I|:=i_1+\cdots +i_{n}$ the order of the differential operator
$X^I$ and  $d(I):=d_1i_1+\cdots +d_ni_{n}$ its degree of\linebreak $\mathbb{G}$-homogeneity. Finally, we denote by $\Delta_\mathbb{G}$ the positive sub-Laplacian associated to $X_1,\dots,X_{m}$
$$
\Delta_\mathbb{G}\coloneqq\sum_{j=1}^m X_j^2.
$$
It is easy to see that
$$
\Delta_\mathbb{G}u =  \divg (\nabla_\mathbb{G} u).
$$
Moreover, $\Delta_\mathbb{G}$ is left-invariant, i.e., for any $x\in\mathbb{G}$, we have
$$\Delta_\mathbb{G} (u\circ \tau_x) = (\Delta_\mathbb{G}u)\circ \tau_x.$$

In parallel, if $E\subset\mathbb{G}$ is a measurable set, a notion of $\mathbb{G}$-perimeter measure
$|\partial E|_\mathbb{G}$ has been introduced in \cite{GN}. We refer to \cite{GN},
\cite{FSSC_houston}, \cite{FSSC_CAG}, \cite{FSSC_step2} for a detailed presentation. For
our needs, we restrict ourselves to recall that if $E$ has locally finite $\mathbb{G}$-perimeter measure
(i.e., it is a $\mathbb{G}$-Caccioppoli set),
then $|\partial E|_\mathbb{G}$ is a Radon measure in $\mathbb{G}$, $\mathbb{G}$-homogeneous of degree $Q-1$. Moreover, the following
representation theorem holds (see  \cite{capdangar}).

\begin{prop}\label{perimetro regolare}
Let $\Omega\subset\mathbb{G}$ be an open set. If $E$ is a $\mathbb{G}$-Caccioppoli set with Euclidean ${\mathbf C}^1$
boundary, then there is an explicit representation of the
$\mathbb{G}$-perimeter in terms of the Euclidean $(N-1)$-dimensional
Hausdorff measure $\mathcal H^{N-1}$
\begin{equation*}
|\partial E|_\mathbb{G}(\Omega)=\int_{\partial
E\cap\Omega}\bigg(\sum_{j=1}^{m_1}\langle
X_j,n\rangle\hspace{0.025cm}_{\mathbb{R}^N}^2\bigg)^{1/2}d{\mathcal {H}}^{N-1},
\end{equation*}
where $n=n(x)$ is the Euclidean unit outward normal to $\partial
E$.
\end{prop}

We also have the subsequent result, which generalizes the classical divergence theorem.
\begin{prop}\label{divergence}
If $E$ is a regular bounded open set with Euclidean ${\mathbf C}^1$
boundary and $\phi$ is a horizontal vector field,
continuously differentiable on $ \overline{E} $, then
$$
\int_{E} \divg \phi\,\hspace{0.075cm}dx = \int_{\partial E} \scal{\phi}{n_\mathbb{G}}\hspace{0.075cm}d |\partial E|_\mathbb{G},
$$
where $n_\mathbb{G}(x)$ is the intrinsic horizontal outward normal to $\partial E$,
given by the (normalized) projection of $n(x)$ on the fiber $H\mathbb{G}_x$ of
the horizontal fiber bundle $H\mathbb{G}$.
\end{prop}

\begin{rem}
The definition of $n_\mathbb{G}$ is well done, since $H\mathbb{G}_x$ is transversal to
the tangent space to $E$ at $x$ for  $|\partial E|_\mathbb{G}$-a.e. $x\in\partial E$
(see \cite{magnani}).
\end{rem}

\medskip

\begin{defn}\label{distanza}{\bf (Carnot-Carath\'eodory
distance)}
An absolutely continuous curve $\gamma:[0,T]\to \mathbb{G}$ is a {\it
sub-unit curve} with respect to $X_1,\dots,X_{m}$ if it is a {\it horizontal curve}, i.e., if there are real measurable
functions $c_1(s),\dots,c_{m}(s)$, $s\in [0,T],$ such that
$$\dot\gamma(s)=\sum\limits_{j=1}^{m}\,c_j(s) X_j(\gamma(s)),
\qquad \text{for a.e.}\;s\in [0,T],$$ and, in addition,
 $$\sum_jc_j^2\le 1.$$

Given $x,y\in\mathbb{G}$,
their Carnot-Carath\'eodory
distance (cc-distance) $d_c(x,y)$ is defined as follows:
$$
d_c(x,y)\coloneqq\inf\left\{T>0:\;\text{there is a sub-unit curve}\;
\gamma\;\text{with}\; \gamma(0)=x,\,\gamma(T)=y\right\}. $$
\end{defn}
The set of sub-unit curves joining $x$ and $y$ is not empty, by
Chow-Rashevsky Theorem (see for instance \cite{BLU}, Theorem 9.1.3). Indeed, by (\ref{stratificazione}), the dimension of the
Lie algebra generated by $X_1,\dots,X_{m}$ is $n.$ Hence, $d_{c}$
is a distance on $\mathbb{G}$ inducing the same topology of the standard
Euclidean distance.

We shall denote by $B_c(x,r)$
 the open  balls associated with $d_c$.   
The cc-distance is well-behaved with respect to left
translations and dilations, that is
\begin{equation}\label{well behaved}
d_c(z\circ x,z\circ y)=d_c(x,y),\quad
d_c(\delta_\lambda(x),\delta_\lambda(y))=\lambda d_c(x,y),
\end{equation}
with $x,y,z\in\mathbb{G}$ and $\lambda>0.$ Using this, we have
 \begin{equation*}\label{ball measures}
|B_c(x,r)| = r^Q |B_c(0,1)| \quad\mbox{and}\quad |\partial B_c(x,r)|_\mathbb{G} (\mathbb{G})
= r^{Q-1} |\partial B_c(0,1)|_\mathbb{G} (\mathbb{G}).
\end{equation*}

	Here and in the sequel, we suppose to work with a real sub-Laplacian $\Delta_{\mathbb{G}} $ on a Carnot group $(\mathbb{G}, \circ,\delta_\lambda)$ with homogeneous dimension $Q\geq 3$. We know, see \cite{BLU}, that there exists a homogeneous norm $d(\cdot)=|\cdot|$ such that
	$$
	\Gamma_{\mathbb{G}} (x,y)=c_Qd^{2-Q}(x^{-1}\circ y),
	$$
	where $\Gamma_{\mathbb{G}}$ is the fundamental solution of $\Delta_\mathbb{G}$. This means that, for every $\varphi\in C_0^{\infty}(\mathbb{R}^N),$ $\Gamma_{\mathbb{G}}:\mathbb{G}\setminus\{y\}\to \mathbb{R}$ satisfies
	\begin{equation*}
	\int_{\mathbb{R}^N}\Gamma_{\mathbb{G}} (x,y)\varphi(x)\hspace{0.075cm}dx=-\varphi(y),
	\end{equation*}
	see Section 5.3 in \cite{BLU} for the details. In order to simplify the notation, if no confusion arises, we may write $\Gamma(0,\xi)\equiv \Gamma(\xi)$ instead of $\Gamma_{\mathbb{G}}(0,\xi)\equiv \Gamma_{\mathbb{G}}(\xi),$ since $\Gamma_{\mathbb{G}}$ is homogeneous. We point out that, in general, the homogeneous norm does not coincide with the Carnot-Cath\'eodory norm $|x|_c\coloneqq d_c(x,0),$ descending from the Carnot-Charath\'eodory distance. Nevertheless, on compact sets they are equivalent, see \cite{Folland_Stein} or \cite{BLU}. In the trivial Euclidean case, it holds $|x|_c=|x|$, possibly up to a multiplicative positive constant.  
	
	Mainly following the notation of Section 5.5 in \cite{BLU}, for every $x_0\in \mathbb{G}$ and $r>0$, we define $D(x_0,r):=\{y\in \mathbb{G}:\hspace{0.2cm}|x_0^{-1}\circ y|<r\}.$ Moreover, for every continuous function $u\in C(\mathbb{G})$ and $r>0,$ we define the following mean value associated with $u$ and $r$ at $x_0$ in the Carnot group $\mathbb{G}$:
	\begin{equation}\label{averageformula}
	\mathcal{M}_r(u)(x_0):=\int_{\partial D(x_0,r)}u(y)\mathcal{K}(x_0^{-1}\circ y)d\sigma(y),
	\end{equation}
	where 
	\begin{equation}\label{def-K}
		\mathcal{K}:=\frac{|\nabla_{\mathbb{G}}\Gamma|^2}{|\nabla \Gamma |}.
	\end{equation}
    We recall a couple of useful examples. In the trivial case of the Euclidean space, it holds, for $x_0,y\in \mathbb{R}^{N},$
    $$\mathcal{K}(y-x_0)=\frac{|\nabla_{\mathbb{R}^N}\Gamma (|y-x_0|)|^2}{|\nabla \Gamma (|y-x_0|)|}=c_N(N-2)|y-x_0|^{1-N}=\frac{r^{1-N}}{N\omega_N },$$ obtaining the classical mean value formula, see for instance \cite{GT}.

In the Heisenberg case $\mathbb{H}^1\equiv \mathbb{R}^3$, taking $X=\frac{\partial}{\partial x}+2y\frac{\partial}{\partial t}$  $Y=\frac{\partial}{\partial y}-2x\frac{\partial}{\partial t}$ as a basis of the horizontal layer and assuming for simplicity $x_0=(0,0,0)$, we obtain, whenever $(x,y,t)\in \mathbb{H}^1,$
\begin{equation*}
\begin{split}\mathcal{K}(x,y,t)=\frac{|\nabla_{\mathbb{H}^1}\Gamma_{\mathbb{H}^1} (x,y,t)|^2}{|\nabla \Gamma_{\mathbb{H}^1} (x,y,t)|}=c|(x,y,t)|_{\mathbb{H}^1}^{-3}\frac{|\nabla_{\mathbb{H}^1}|(x,y,t)|_{\mathbb{H}^1}|^2}{|\nabla |(x,y,t)|_{\mathbb{H}^1}|},
\end{split}
\end{equation*}
where $\Gamma_{\mathbb{H}^1}(x,y,t)=c|(x,y,t)|_{\mathbb{H}^1}^{-2}$, $|(x,y,t)|_{\mathbb{H}^1}=((x^2+y^2)^2+t^2)^{\frac{1}{4}}$, $$\nabla_{\mathbb{H}^1}|(x,y,t)|_{\mathbb{H}^1}=X|(x,y,t)|_{\mathbb{H}^1}X+ Y|(x,y,t)|_{\mathbb{H}^1}Y,$$ and 
$$|\nabla_{\mathbb{H}^1}|(x,y,t)|_{\mathbb{H}^1}|=\sqrt{(X|(x,y,t)|_{\mathbb{H}^1})^2+ (Y|(x,y,t)|_{\mathbb{H}^1})^2}.$$

	We are in position now to recall the main result we use in our proof, which characterizes the harmonic and subharmonic functions in Carnot groups, see \cite{BLU} for the details.
\begin{prop}[\cite{BLU}]\label{Bonfi}
Let $u\in C(\Omega)$, where $\Omega\subset \mathbb{G}$ is an open connected set. Then, $u\in C^{\infty}(\Omega),$ $\Delta_{\mathbb{G}}u=0$ if and only if for every $\overline{D}(x,r)\subset \Omega$
$$
u(x)=\mathcal{M}_r(u)(x).
$$  
In addition, a function $u$ is subharmonic if an only if, for any $x\in \Omega,$ $r\to \mathcal{M}_r(u)(x)$ is monotone increasing for $0<r\le R(x)$ and 
$$
u(x)=\lim_{r\to 0^+}\mathcal{M}_r(u)(x).
$$
\end{prop}	
We conclude this introductory part by recalling that, like in the classical case, the following result even holds in the noncommutative framework given by Carnot groups.
\begin{prop}\label{propoeasy}
Let $\Omega\subset \mathbb{G}$ be an open set. It holds: 
\begin{itemize}
\item[(i)] If $u,v$ are two continuous subharmonic functions, then $\max\{u,v\}$ is a subharmonic function.
\item[(ii)] If $u,v$ are two continuous superharmonic functions, then  $\min\{u,v\}$ is a  superharmonic function.
\end{itemize}
\end{prop}
\begin{proof}
Let us prove the case (i). Since for every $D(x_0,r)\subset \Omega$
$$
u(x_0)\leq \mathcal{M}_r(u)(x_0),\quad v(x_0)\leq \mathcal{M}_r(v)(x_0),
$$
then
\begin{equation*}\begin{split}
&\mathcal{M}_r(\max\{u,v\})(x_0)=\int_{\partial D(x_0,r)}\max\{u,v\}(y)\mathcal{K}(x_0^{-1}\circ y)d\sigma(y)\\
&\geq \int_{\partial D(x_0,r)}u(y)\mathcal{K}(x_0^{-1}\circ y)d\sigma(y)\geq u(x_0)
\end{split}
\end{equation*}
and, arguing in the same way with $v$ instead of $u,$ 
\begin{equation*}
\mathcal{M}_r(\max\{u,v\})(x_0)\geq v(x_0).
\end{equation*}
As a consequence, we have
\begin{equation*}
\mathcal{M}_r(\max\{u,v\})(x_0)\geq \max\{u,v\}(x_0).
\end{equation*}
The case (ii) follows by a similar argument.
\end{proof}
\begin{rem}
Let u be a continuous harmonic function in a Carnot group $\mathbb{G}.$ Then, $u^+:=\max\{u,0\}$ and $u^-:=\max\{-u,0\}$ are both subharmonic. This is an easy consequence of previous Proposition \ref{propoeasy}.
\end{rem}
\section{Proofs of main results}\label{treF}
In this section, we state our main result in its general form. It is a representation formula of the function $J_{u,x_0}^{\mathbb{G}}$ that guarantees some simple consequences thanks to correct interpretations of the players inside it, mainly the average formula. In other words, we present Theorem \ref{maintheorem} in the following more detailed way.
\begin{thm}\label{maintheoremgeneral}
		Let $\Omega\subset \mathbb{G}$ be an open set. Let $u\in C(\Omega)$ such that \linebreak $|\nabla_{\mathbb{G}} u|^2\in L^1_{\mathrm{loc}}(\Omega)$.		
		For every $x_0\in \partial \{u>0\}\cap \Omega$, 
		it holds
		\begin{equation*}\label{moneqH5starbisstat}\begin{split}
		J_{u,x_0}^{\mathbb{G}}(r)
		&=\frac{1}{Q-2}\int_{0}^1t\mathcal{M}_{tr}\left(\left(\frac{|\nabla_{\mathbb{G}} u^+(x_0\circ\xi)|}{|\nabla_{\mathbb{G}}|\xi||}\right)^2\right)(0)\hspace{0.075cm}dt\\
		&\cdot\frac{1}{Q-2}\int_0^1t\mathcal{M}_{tr}\left(\left(\frac{|\nabla_{\mathbb{G}} u^-(x_0\circ\xi)|}{|\nabla_{\mathbb{G}}|\xi||}\right)^2\right)(0)\hspace{0.075cm}dt.
\end{split}
	\end{equation*} 
	Consequently, if $$t\to \mathcal{M}_{t}\left(\left(\frac{|\nabla_{\mathbb{G}} u^\pm(x_0\circ\xi)|}{|\nabla_{\mathbb{G}}|\xi||}\right)^2\right)(0)$$  are monotone increasing almost everywhere, then the function $J_{u, x_0}^{\mathbb{G}}$ is monotone increasing in a right neighborhood of $0$.
	
	In particular, for every $x_0\in \Omega,$ we have
	\begin{equation*}\label{moneqH5starbisstat}\begin{split}
		 \frac{1}{r^2}\int_{D(x_0,r)}|\nabla_{\mathbb{G}} u|^2\Gamma(x_0,\xi)\hspace{0.075cm}d\xi=\frac{1}{Q-2}\int_{0}^1t\mathcal{M}_{tr}\left(\left(\frac{|\nabla_{\mathbb{G}} u(x_0\circ\xi)|}{|\nabla_{\mathbb{G}}|\xi||}\right)^2\right)(0)\hspace{0.075cm}dt.
\end{split}
	\end{equation*} 
	Thus, if $$t\to \mathcal{M}_{t}\left(\left(\frac{|\nabla_{\mathbb{G}} u(x_0\circ\xi)|}{|\nabla_{\mathbb{G}}|\xi||}\right)^2\right)(0)$$  is monotone increasing a.e., then
	$$
	r\to \frac{1}{r^2}\int_{D(x_0,r)}|\nabla_{\mathbb{G}} u|^2\Gamma(x_0,\xi)\hspace{0.075cm}d\xi
	$$
	is monotone increasing in a right neighborhood of $0$.
		\end{thm}
%
\begin{proof}[Proof of Theorem \ref{maintheoremgeneral} and Theorem \ref{maintheorem}]
After a left translation, we may rewrite \eqref{moneqH01} as
\begin{equation*}\label{moneqH}
		J_{u}^{\mathbb{G}}(r)= \frac{1}{r^4}\int_{D(0,r)}|\nabla_{\mathbb{G}} u^+|^2\hspace{0.075cm}\Gamma (0,\xi)\hspace{0.075cm}d\xi\int_{D(0,r)}|\nabla_{\mathbb{G}} u^-|^2\hspace{0.075cm}\Gamma (0,\xi)\hspace{0.075cm}d\xi,
	\end{equation*}
	that is we may assume $x_0=0$ and $u=u(x_0\hspace{0.05cm}\circ\hspace{0.05cm}\cdot)$, denoting by $J_{u,0}^{\mathbb{G}}(r)\equiv J_{u}^{\mathbb{G}}(r)$.\\
	We remark that, defining $u_{r}(x):=\frac{u(\delta_r( x))}{r}$, it holds
\begin{equation}\label{moneqH1}\begin{split}
		&J_{u}^{\mathbb{G}}(r)=\\
		&=\frac{1}{r^4}\int_{D(0,1)}|\nabla_{\mathbb{G}} u^+(\delta_r\xi)|^2\Gamma(0,\delta_r\xi)r^{Q}\hspace{0.075cm}d\xi\int_{D(0,1)}|\nabla_{\mathbb{G}} 
		u^-(\delta_r\xi)|^2\Gamma(0,\delta_r \xi)\hspace{0.05cm}r^Q\hspace{0.075cm}d\xi\\
		&= \frac{1}{r^4}\int_{D(0,1)}|\nabla_{\mathbb{G}} u_r^+(\xi)|^2\hspace{0.075cm}\Gamma(0,\xi)r^{2}\hspace{0.075cm}d\xi\int_{D(0,1)}|\nabla_{\mathbb{G}} u_r^-(\xi)|^2\hspace{0.075cm}\Gamma(0,\xi)\hspace{0.05cm}r^2\hspace{0.075cm}d\xi\\
		&=J_{u_r}^{\mathbb{G}}(1).
		\end{split}
	\end{equation}
	In addition, for every $\lambda>0$, the function $u_{\lambda}$ satisfies $\Delta_{\mathbb{G}}u_\lambda=\lambda (\Delta_{\mathbb{G}}u)(\delta_\lambda(x))$  and $\nabla_{\mathbb{G}} u_{\lambda}(x)=(\nabla_{\mathbb{G}} u)(\delta_\lambda x)$ in $\delta_{\lambda^{-1}}(\Omega)$.\\
	Let us study now the behavior of $ J_{u_r}^{\mathbb{G}}(1),$ with
\begin{equation*}\label{moneqH2}\begin{split}
		J_{u_r}^{\mathbb{G}}(1)&=\int_{D(0,1)}|\nabla_{\mathbb{G}} u_r^+|^2\hspace{0.075cm}\Gamma(0,\xi)\hspace{0.075cm}d\xi\int_{D(0,1)}|\nabla_{\mathbb{G}} u_r^-|^2\hspace{0.075cm}\Gamma(0,\xi)\hspace{0.075cm}d\xi.
		\end{split}
	\end{equation*}
	Recalling the coarea formula, we first obtain
	\begin{equation*}\label{moneqH3}\begin{split}
		J_{u_r}^{\mathbb{G}}(1)&=\int_{0}^1\left(\int_{\partial D(0,t)}|\nabla_{\mathbb{G}} u_r^+|^2\frac{\Gamma(0,\xi)}{|\nabla_{\mathbb{G}}|\xi||}d\mathcal{H}^{Q-1}_{\mathbb{G}}(\xi)\right)dt\\
		&\cdot\int_0^1\left(\int_{\partial D(0,t)}|\nabla_{\mathbb{G}} u_r^-|^2\frac{\Gamma(0,\xi)}{|\nabla_{\mathbb{G}}|\xi||} \hspace{0.05cm}d\mathcal{H}^{Q-1}_{\mathbb{G}}(\xi)\right)dt\\
		&=c_Q\int_{0}^1t^{2-Q}\left(\int_{\partial D(0,t)}|\nabla_{\mathbb{G}} u_r^+|^2\frac{1}{|\nabla_{\mathbb{G}}|\xi||}d\mathcal{H}^{Q-1}_{\mathbb{G}}(\xi)\right)dt\\
		&\cdot c_Q\int_0^1t^{2-Q}\left(\int_{\partial D(0,t)}|\nabla_{\mathbb{G}} u_r^-|^2\frac{1}{|\nabla_{\mathbb{G}}|\xi||} \hspace{0.05cm}d\mathcal{H}^{Q-1}_{\mathbb{G}}(\xi)\right)dt,
		\end{split}
	\end{equation*}
	in other words
	\begin{equation*}\label{moneqH4}\begin{split}
		J_{u_r}^{\mathbb{G}}(1)&=c_Q\int_{0}^1\frac{1}{t^{Q-2}}\left(\int_{\partial D(0,t)}|\nabla_{\mathbb{G}} u_r^+|^2\frac{1}{|\nabla_{\mathbb{G}}|\xi||}d\mathcal{H}^{Q-1}_{\mathbb{G}}(\xi)\right)dt\\
		&\cdot c_Q\int_0^1\frac{1}{t^{Q-2}}\left(\int_{\partial D(0,t)}|\nabla_{\mathbb{G}} u_r^-|^2\frac{1}{|\nabla_{\mathbb{G}}|\xi||}\hspace{0.05cm}d\mathcal{H}^{Q-1}_{\mathbb{G}}(\xi)\right)dt.
		\end{split}
	\end{equation*} 
By definition of $u_r,$ this implies
\begin{equation*}\begin{split}
		&J_{u_r}^{\mathbb{G}}(1)=c_Q\int_{0}^1\frac{1}{t^{Q-2}}\left(\int_{\partial D(0,t)}|\nabla_{\mathbb{G}} u^+(\delta_{r}(\xi))|^2\frac{1}{|\nabla_{\mathbb{G}}|\xi||}d\mathcal{H}^{Q-1}_{\mathbb{G}}(\xi)\right)dt\\
		&\cdot c_Q\int_0^1\frac{1}{t^{Q-2}}\left(\int_{\partial D(0,t)}|\nabla_{\mathbb{G}} u^-(\delta_{r}(\xi))|^2\frac{1}{|\nabla_{\mathbb{G}}|\xi||}\hspace{0.05cm}d\mathcal{H}^{Q-1}_{\mathbb{G}}(\xi)\right)dt.
		\end{split}
		\end{equation*}
		Applying a change of variables in the surface integral, we then achieve
		\begin{equation*}\begin{split}\label{moneqH42}
		&J_{u_r}^{\mathbb{G}}(1)=c_Q\int_{0}^1\frac{t}{(tr)^{Q-1}}\left(\int_{\partial D(0,tr)}|\nabla_{\mathbb{G}} u^+|^2\frac{1}{|\nabla_{\mathbb{G}}|z||(\delta_{r^{-1}}(\xi))}d\mathcal{H}^{Q-1}_{\mathbb{G}}(\xi)\right)dt\\
		&\cdot c_Q\int_0^1\frac{t}{(tr)^{Q-1}}\left(\int_{\partial D(0,tr)}|\nabla_{\mathbb{G}} u^-(\xi)|^2\frac{1}{|\nabla_{\mathbb{G}}|z||(\delta_{r^{-1}}(\xi))}\hspace{0.05cm}d\mathcal{H}^{Q-1}_{\mathbb{G}}(\xi)\right)dt.
		\end{split}
	\end{equation*} 
	On the other hand, by homogeneity, $\nabla_{\mathbb{G}}|z||(\delta_{r^{-1}}\xi)= \nabla_{\mathbb{G}}|\xi|.$	Hence, we get
	\begin{equation}\label{moneqH5}\begin{split}
		J_{u_r}^{\mathbb{G}}(1)&=c_Q\int_{0}^1\frac{t}{(tr)^{Q-1}}\left(\int_{\partial D(0,tr)}|\nabla_{\mathbb{G}} u^+|^2\frac{1}{|\nabla_{\mathbb{G}}|\xi||}d\mathcal{H}^{Q-1}_{\mathbb{G}}(\xi)\right)dt\\
		&\cdot c_Q\int_0^1\frac{t}{(tr)^{Q-1}}\left(\int_{\partial D(0,tr)}|\nabla_{\mathbb{G}} u^-|^2\frac{1}{|\nabla_{\mathbb{G}}|\xi||}\hspace{0.05cm}d\mathcal{H}^{Q-1}_{\mathbb{G}}(\xi)\right)dt.
		\end{split}
	\end{equation} 
	We focus now on the equality
	\begin{equation}\label{moneqHG1star}\begin{split}
\frac{1}{|\nabla_{\mathbb{G}}|\xi||}\hspace{0.05cm}d\mathcal{H}^{Q-1}_{\mathbb{G}}(\xi)=\frac{1}{|\nabla|\xi||}\hspace{0.05cm}d\mathcal{H}^{N-1}_{\mathbb{R}^N}(\xi).
	\end{split}
	\end{equation}		
	We first notice that 
\begin{equation*}\label{moneqHG1}\begin{split}
\nabla_{\mathbb{G}}|\xi|&= \nabla_{\mathbb{G}}((c_Q^{-1}\Gamma (0,\xi))^{-\frac{1}{Q-2}})=-\frac{c_Q^{\frac{1}{Q-2}}}{Q-2}\Gamma(0,\xi)^{-\frac{1}{Q-2}-1}\nabla_{\mathbb{G}}\Gamma(0,\xi)\\
&=-\frac{c_Q^{-1}}{Q-2}|\xi|^{Q-1}\nabla_{\mathbb{G}}\Gamma(0,\xi),
\end{split}
	\end{equation*}	
	so that 
	\begin{equation}\label{moneqHG18}\begin{split}
\frac{|\nabla_{\mathbb{G}}|\xi||}{|\nabla_{\mathbb{G}}\Gamma(0,\xi)|}&=\frac{c_Q^{-1}}{Q-2}|\xi|^{Q-1}.
\end{split}
	\end{equation}	
	Repeating the same argument with $\nabla$ instead of $\nabla_{\mathbb{G}},$ it holds
	\begin{equation*}\label{moneqHG22}\begin{split}
\nabla |\xi|=-\frac{c_Q^{-1}}{Q-2}|\xi|^{Q-1}\nabla\Gamma(0,\xi).
\end{split}
	\end{equation*}	
	Therefore, using this and \eqref{moneqHG18}, it follows, from \eqref{moneqHG1star},
	\begin{equation}\label{moneqHG1star2}\begin{split}
&\frac{1}{|\nabla_{\mathbb{G}}|\xi||}\hspace{0.05cm}d\mathcal{H}^{Q-1}_{\mathbb{G}}(\xi)=\frac{|\nabla_{\mathbb{G}}\Gamma(0,\xi)|^2}{|\nabla_{\mathbb{G}}\Gamma(0,\xi)|^2|\nabla|\xi||}\hspace{0.05cm}d\mathcal{H}^{N-1}_{\mathbb{R}^N}(\xi)\\
&=\frac{Q-2}{c_Q^{-1}}\frac{|\nabla_{\mathbb{G}}\Gamma(0,\xi)|^2}{|\nabla_{\mathbb{G}}\Gamma(0,\xi)|^2|\xi|^{Q-1}|\nabla\Gamma(0,\xi)|}\hspace{0.05cm}d\mathcal{H}^{N-1}_{\mathbb{R}^N}(\xi)\\
&= \frac{Q-2}{c_Q^{-1}}\frac{\mathcal{K}}{|\nabla_{\mathbb{G}}\Gamma(0,\xi)|^2|\xi|^{Q-1}}\hspace{0.05cm}d\mathcal{H}^{N-1}_{\mathbb{R}^N}(\xi).
	\end{split}
	\end{equation}	
	On the other hand, since $$
	\Gamma_{\mathbb{G}} (x,y)=c_Qd^{2-Q}(x^{-1}\circ y)=c_Q|x^{-1}\circ y|^{2-Q},
	$$
	we deduce that
	$$
	\nabla_{\mathbb{G}}\Gamma (0,\xi)=c_Q(2-Q)|\xi|^{1-Q}\nabla_{\mathbb{G}}|\xi|.
	$$ 
	Thus, we have, by \eqref{moneqHG1star2},
	\begin{equation*}\label{moneqHG1star22}\begin{split}
&\frac{1}{|\nabla_{\mathbb{G}}|\xi||}\hspace{0.05cm}d\mathcal{H}^{Q-1}_{\mathbb{G}}(\xi)= \frac{1}{(Q-2)c_Q}\frac{|\xi|^{Q-1}}{|\nabla_{\mathbb{G}}|\xi||^2}\hspace{0.05cm} \mathcal{K} d\mathcal{H}^{N-1}_{\mathbb{R}^N}(\xi).
	\end{split}
	\end{equation*}
	Now, recalling \eqref{moneqH5} and keeping in mind this equality, we obtain 
	\begin{equation*}\label{moneqH5star}\begin{split}
		J_{u_r}^{\mathbb{G}}(1)&=c_Q\int_{0}^1\frac{t}{(tr)^{Q-1}}\left(\int_{\partial D(0,tr)}\left(\frac{|\nabla_{\mathbb{G}} u^+|}{|\nabla_{\mathbb{G}}|\xi||}\right)^2 \frac{|\xi|^{Q-1}\mathcal{K}}{(Q-2)c_Q}\hspace{0.05cm}d\mathcal{H}^{N-1}_{\mathbb{R}^N}(\xi)\right)dt\\
		&\cdot c_Q\int_{0}^1\frac{t}{(tr)^{Q-1}}\left(\int_{\partial D(0,tr)}\left(\frac{|\nabla_{\mathbb{G}} u^-|}{|\nabla_{\mathbb{G}}|\xi||}\right)^2 \frac{|\xi|^{Q-1}\mathcal{K}}{(Q-2)c_Q}\hspace{0.05cm}d\mathcal{H}^{N-1}_{\mathbb{R}^N}(\xi)\right)dt\\
		&=\frac{1}{Q-2}\int_{0}^1t\left(\int_{\partial D(0,tr)}\left(\frac{|\nabla_{\mathbb{G}} u^+|}{|\nabla_{\mathbb{G}}|\xi||}\right)^2\mathcal{K}\hspace{0.05cm}d\mathcal{H}^{N-1}_{\mathbb{R}^N}(\xi)\right)dt\\
		&\cdot\frac{1}{Q-2}\int_0^1t\left(\int_{\partial D(0,tr)}\left(\frac{|\nabla_{\mathbb{G}} u^-|}{|\nabla_{\mathbb{G}}|\xi||}\right)^2\mathcal{K}\hspace{0.05cm}d\mathcal{H}^{N-1}_{\mathbb{R}^N}(\xi)\right)dt.
		\end{split}
	\end{equation*} 
	This, according to \eqref{averageformula}, yields
	\begin{equation*}\label{moneqH5starbis}\begin{split}
		J_{u_r}^{\mathbb{G}}(1)
		&=\frac{1}{Q-2}\int_{0}^1t\mathcal{M}_{tr}\left(\left(\frac{|\nabla_{\mathbb{G}} u^+|}{|\nabla_{\mathbb{G}}|\xi||}\right)^2\right)(0)dt\\
		&\cdot\frac{1}{Q-2}\int_0^1t\mathcal{M}_{tr}\left(\left(\frac{|\nabla_{\mathbb{G}} u^-|}{|\nabla_{\mathbb{G}}|\xi||}\right)^2\right)(0)dt.
\end{split}
	\end{equation*}
    Using \eqref{moneqH1}, we finally achieve the desired expression of $J_{u}^{\mathbb{G}}.$\\
	Since by hypothesis $s\to \mathcal{M}_{s}\left(\left(\frac{|\nabla_{\mathbb{G}} u^{\pm}(\xi)|}{|\nabla_{\mathbb{G}}|\xi||}\right)^2\right)(0)$ are monotone increasing a.e., 
	for every $r_1\leq r_2,$ we have, for almost every $t\in [0,1],$ 
	$$
	t\mathcal{M}_{tr_1}\left(\left(\frac{|\nabla_{\mathbb{G}} u^{\pm}(\xi)|}{|\nabla_{\mathbb{G}}|\xi||}\right)^2\right)(0)\leq t\mathcal{M}_{tr_{2}}\left(\left(\frac{|\nabla_{\mathbb{G}} u^-(\xi)|}{|\nabla_{\mathbb{G}}|\xi||}\right)^2\right)(0).
	$$
	Hence, we get
	$$
	J_{u}^{\mathbb{G}}(r_1)=J_{u_{r_1}}^{\mathbb{G}}(1)\leq J_{u_{r_2}}^{\mathbb{G}}(1)=J_{u}^{\mathbb{G}}(r_2),
	$$
	concluding the proof. The case in which we have only one function is a straightforward consequence of the previous arguments. 
	\end{proof}

	\begin{rem}\label{rem-ACF}
We remark that, recalling Proposition \ref{Bonfi}, in case $$\xi\to \left(\frac{|\nabla_{\mathbb{G}} u^{\pm}(x_0\circ\xi)|}{|\nabla_{\mathbb{G}}|\xi||}\right)^2$$ were
subharmonic,
we would immediately obtain the monotone increasing behavior of $s\to \mathcal{M}_{s}\left(\left(\frac{|\nabla_{\mathbb{G}} u^{\pm}(x_0\circ\xi)|}{|\nabla_{\mathbb{G}}|\xi||}\right)^2\right)(0),$ and, as a consequence, the Alt-Caffarelli-Friedman monotonicity formula whenever $x_0\in \partial\{u>0\}\cap \Omega$.
   	\end{rem}
	Even in the particular and classical case in which $\mathbb{G}=\mathbb{R}^N$ is a group of step one, that is in the Euclidean framework, we obtain $\nabla_{\mathbb{G}}|\xi|=\nabla |\xi|=\frac{\xi}{|\xi|}$. Hence, it holds, in view of \eqref{moneqH1},
	\begin{equation*}\label{moneqH7}\begin{split}
		J_{u}^{\mathbb{\mathbb{R}^N}}(r)=J_{u_r}^{\mathbb{\mathbb{R}^N}}(1)&=\frac{1}{(N-2)^2}\int_{0}^1t\mathcal{M}_{tr}(|\nabla u^+|^2)(0)\hspace{0.05cm}dt \int_{0}^1t\mathcal{M}_{tr}(|\nabla u^-|^2)(0)\hspace{0.05cm}dt.
		\end{split}
	\end{equation*} 
	Unfortunately, in general, $|\nabla u^+|^2$ and $|\nabla u^-|^2$ are not subharmonic, so we can not achieve the Alt-Caffarelli-Friedman monotonicity formula exploiting Remark \ref{rem-ACF}. However, in the simpler case described in Corollary \ref{yes}, this argument applies.
	
\begin{proof}[Proof of Corollary \ref{yes}] Since $u$ is harmonic, for every $i=1\dots,n$ $\frac{\partial u}{\partial x_i}$ is harmonic as well. Thus, by a straightforward computation, $(\frac{\partial u}{\partial x_i})^2$ is subharmonic, which obviously implies that $|\nabla u|^2$ is subharmonic. Hence, recalling that in this case
$$\left(\frac{|\nabla u(x_0+\xi)|}{|\nabla |\xi||}\right)^2=|\nabla u(x_0+\xi)|^2,$$ 
we obtain, from Theorem \ref{maintheoremgeneral},
\begin{equation*}\label{moneqH7}\begin{split}
	\frac{c_N}{r^2}\int_{B(0,r)}\frac{|\nabla u(x_0+\xi)|^2}{|\xi|_{\mathbb{R}^N}^{N-2}}\hspace{0.075cm}d\xi		=\frac{1}{(N-2)^2}\int_{0}^1t\mathcal{M}_{tr}(|\nabla u(x_0+\cdot)|^2)(0)\hspace{0.05cm}dt.
		\end{split}
		\end{equation*}
		Since  $|\nabla u(x_0+\cdot)|^2$ is subharmonic, we have that 
		$$
		r\to\mathcal{M}_{r}(|\nabla u(x_0+\cdot)|^2
		$$
		is monotone increasing, proving the monotone increasing behavior of 
		$$
		r\to \frac{c_N}{r^2}\int_{B(0,r)}\frac{|\nabla u(x_0+\xi)|^2}{|\xi|_{\mathbb{R}^N}^{N-2}}\hspace{0.075cm}d\xi	.
		$$
\end{proof}
	\begin{proof}[Proof of Corollary \ref{maintheorem2}]
	In the simplest noncommutative case of the Heisenberg group $\mathbb{H}^1$, in the recent paper \cite{FeFo2}, the authors have proved that the Alt-Caffarelli-Friedman monotonicity formula determined by $J_{u}^{\mathbb{H}^1}$ fails. In that paper, we have produced a counterexample to the monotone increasing behavior of $J_{u}^{\mathbb{H}^1}.$ To get it, we have followed a different idea associated with the lack of orthogonality of harmonic polynomials in the Heisenberg group.
	
More precisely, we have first proved that, given the polynomial \begin{equation}\label{polyn-count}
	u=x-3yt-2x^3,
\end{equation} 
	it turns out to be harmonic (a simple straightforward computation) and we find that
	$$
	\tilde{J}_{u}^{\mathbb{H}^1}(r)\coloneqq \frac{1}{r^2}\int_{D(0,r)}|\nabla_{\mathbb{G}} u|^2\hspace{0.075cm}\Gamma (0,\xi)\hspace{0.075cm}d\xi
	$$
	is monotone decreasing in a right neighborhood of $r=0.$
	
	The straightforward proof of the decreasing monotonicity property of $\tilde J_{u}^{\mathbb{H}^1}$ appears a little more delicate in this case, due to the structure of $\mathcal{K}$. Thus, we first argue by introducing a computation on the solid integral, see \cite{FeFo2}. Precisely, let us rewrite $\tilde J_{u}^{\mathbb{H}^1}$ as 
	\begin{equation*}\label{moneqHcom}\begin{split}
		&\tilde{J}_{u}^{\mathbb{H}^1}(r)=\frac{1}{r^2}\int_{D(0,r)\cap\{u>0\}}|\nabla_{\mathbb{G}} u|^2\hspace{0.075cm}\Gamma (0,\xi)\hspace{0.075cm}d\xi\\
		&+\frac{1}{r^2}\int_{D(0,r)\cap\{u<0\}}|\nabla_{\mathbb{G}} u|^2\hspace{0.075cm}\Gamma (0,\xi)\hspace{0.075cm}d\xi.
		\end{split}
	\end{equation*}
Hence, performing the change of variables $(x,y,t)\to(-x,-y,t)$  in the second integral, see \cite{FeFo2}, we obtain
	\begin{equation*}\label{moneqHcom2}\begin{split}
		&\tilde{J}_{u}^{\mathbb{H}^1}(r)=\frac{1}{r^2}\int_{D(0,r)\cap\{u>0\}}|\nabla_{\mathbb{G}} u|^2\hspace{0.075cm}\Gamma (0,\xi)\hspace{0.075cm}d\xi\\
		&+\frac{1}{r^2}\int_{D(0,r)\cap\{u>0\}}|\nabla_{\mathbb{G}} u|^2\hspace{0.075cm}\Gamma (0,\xi)\hspace{0.075cm}d\xi,\\
		\end{split}
	\end{equation*}
	because $u(-x,-y,t)=-u(x,y,t)$ by \eqref{polyn-count}. This implies $$\tilde{J}_{u}^{\mathbb{H}^1}(r)=\frac{2}{r^2}\int_{D(0,r)}|\nabla_{\mathbb{G}} u^+|^2\hspace{0.075cm}\Gamma (0,\xi)\hspace{0.075cm}d\xi,$$
	and with the same argument we also have
	\begin{equation*}\label{moneqHcom23}\begin{split}
	\tilde{J}_{u}^{\mathbb{H}^1}(r)= \frac{2}{r^2}\int_{D(0,r)}|\nabla_{\mathbb{G}} u^-|^2\hspace{0.075cm}\Gamma (0,\xi)\hspace{0.075cm}d\xi.
	\end{split}
	\end{equation*}
	Putting together these equalities, it then holds, with the choice of $u$ as in \eqref{polyn-count},
	\begin{equation}\label{moneqHcom24}\begin{split}
		&J_{u}^{\mathbb{H}^1}(r)\coloneqq \frac{1}{4}(\tilde{J}_{u}^{\mathbb{H}^1}(r))^2.
		\end{split}
	\end{equation}
Let us prove that, again with the choice of $u$ as in \eqref{polyn-count}, $\tilde{J}_{u}^{\mathbb{H}^1}$ is monotone decreasing. First, see \cite{FeFo2}, $\tilde{J}_{u}^{\mathbb{H}^1}$ can be rewritten as
	 \begin{equation*}\label{moneqHcom242}\begin{split}
		\tilde{J}_{u}^{\mathbb{H}^1}(r)=\frac{1}{r^2}\int_{D(0,r)}((1 -6x^2-6y^2)^2+9(-t+2xy)^2)\Gamma (0,\xi)\hspace{0.075cm}d\xi.
		\end{split}
	\end{equation*}
	 This yields
	 \begin{equation*}\label{moneqHcom25}\begin{split}
		&\tilde{J}_{u}^{\mathbb{H}^1}(r)=\frac{r^4}{r^2}\int_{D(0,1)}((1 -6r^2(x^2+y^2))^2+9r^4(-t+2xy)^2)r^{-2}\Gamma (0,\xi)\hspace{0.075cm}d\xi\\
		&=\int_{D(0,1)}((1 -6r^2(x^2+y^2))^2+9r^4(-t+2xy)^2)\Gamma (0,\xi)\hspace{0.075cm}d\xi\\
		&=a_0-a_1r^2+a_2r^4,
		\end{split}
	\end{equation*}
	where
	$$
	a_0=\int_{D(0,1)}\Gamma (0,\xi)\hspace{0.075cm}d\xi,\quad a_1=12\int_{D(0,1)}(x^2+y^2)\Gamma (0,\xi)\hspace{0.075cm}d\xi,
	$$
	and
	$$
	a_2=\int_{D(0,1)}(36(x^2+y^2)^2+9(-t+2xy)^2)\Gamma (0,\xi)\hspace{0.075cm}d\xi.
	$$
Hence, $r\to \tilde{J}_{u}^{\mathbb{H}^1}(r)$ is indeed decreasing in a right neighborhood of $0$.\\
Since $\tilde{J}_{u}^{\mathbb{H}^1}(r)$ is positive, by \eqref{moneqHcom24} it then holds that $J_{u}^{\mathbb{H}^1}(r)$ is decreasing.

Let us provide now a different proof for the decreasing behavior of $\tilde{J}_{u}^{\mathbb{H}^1}(r).$ First, applying the same symmetry argument on $u,$ we may reduce ourselves to study the monotonicity behavior, in a right neighborhood of $0,$ of 
\begin{equation}\label{moneqH5starbisx}\begin{split}
		\tilde{J}_{u}^{\mathbb{H}^1}(r)
		&=c_{\mathbb{H}^1}\int_{0}^1s\mathcal{M}_{sr}\left(\left(\frac{|\nabla_{\mathbb{H}^1} u(\xi)|}{|\nabla_{\mathbb{H}^1}|\xi||}\right)^2\right)(0)ds.
		\end{split}
	\end{equation} 
In particular, we focus on
$$
r\to s\mathcal{M}_{sr}\left(\left(\frac{|\nabla_{\mathbb{H}^1} u(\xi)|}{|\nabla_{\mathbb{H}^1}|\xi||}\right)^2\right)(0),
$$
remarking that
\begin{equation*}\begin{split}
		&\left(\mathcal{M}_{sr}\left(\frac{|\nabla_{\mathbb{H}^1} u(\xi)|}{|\nabla_{\mathbb{H}^1}|\xi||}\right)^2\right)(0)=\\
		&=\int_{\partial D(0,sr)}\frac{(1 -6(x^2+y^2))^2+9(-t+2xy)^2}{|\nabla_{\mathbb{H}^1}|\xi||^2}\mathcal{K}(\xi)d\mathcal{H}^2(\xi)\\
		&=\int_{\partial D(0,1)}\frac{(1 -6(sr)^2(x^2+y^2))^2+9(sr)^4(-t+2xy)^2}{|\nabla_{\mathbb{H}^1}|\xi||^2}\mathcal{K}(\xi)d\mathcal{H}^2(\xi)\\
    	\end{split}
\end{equation*}
\begin{equation}\label{moneqH5starbisso}\begin{split}
		&=\mathcal{M}_{1}\left(\frac{1}{|\nabla_{\mathbb{H}^1}|\xi||^2}\right)(0)-12(sr)^2\mathcal{M}_{1}\left(\frac{x^2+y^2}{|\nabla_{\mathbb{H}^1}|\xi||^2}\right)(0)\\
		&+(sr)^4\mathcal{M}_{1}\left(\frac{36(x^2+y^2)^2+9(-t+2xy)^2}{|\nabla_{\mathbb{H}^1}|\xi||^2}\right)(0). 
		\end{split}
	\end{equation} 
The polynomial in $sr$ in \eqref{moneqH5starbisso} is decreasing in a right neighborhood of $0.$ 
Consequently, it holds, for every $0\leq r_1\leq r_2,$
\begin{equation*}\label{moneqH5starbissoc}\begin{split}
		&\mathcal{M}_{sr_1}\left(\left(\frac{|\nabla_{\mathbb{H}^1} u(\xi)|}{|\nabla_{\mathbb{H}^1}|\xi||}\right)^2\right)(0)\geq \mathcal{M}_{sr_2}\left(\left(\frac{|\nabla_{\mathbb{H}^1} u(\xi)|}{|\nabla_{\mathbb{H}^1}|\xi||}\right)^2\right)(0),
		\end{split}
	\end{equation*} 
which implies, for every $s\in [0,1],$
\begin{equation*}\label{moneqH5starbissod}\begin{split}
		&s\mathcal{M}_{sr_1}\left(\left(\frac{|\nabla_{\mathbb{H}^1} u(\xi)|}{|\nabla_{\mathbb{H}^1}|\xi||}\right)^2\right)(0)\geq s\mathcal{M}_{sr_2}\left(\left(\frac{|\nabla_{\mathbb{H}^1} u(\xi)|}{|\nabla_{\mathbb{H}^1}|\xi||}\right)^2\right)(0).
		\end{split}
	\end{equation*} 
As a byproduct, recalling \eqref{moneqH5starbisx}, it follows
$$
\tilde{J}_{u}^{\mathbb{H}^1}(r_1)\geq \tilde{J}_{u}^{\mathbb{H}^1}(r_2),
$$
proving the decreasing monotonicity of $J_{u}^{\mathbb{H}^1},$ exactly arguing as in the solid case.

Now, we can generalize this argument finding several functions with this decreasing property. Indeed, for every $c_1,c_2\in \mathbb{R},$ the function $$u(x,y,t)=c_1x+c_2y+3t(c_2x-c_1y)-2(c_1x^3+c_2y^3)$$ is harmonic in $\mathbb{H}^1$ and $u(-x,-y,t)=-u(x,y,t)$ holds. Thus, applying the same argument, we achieve \eqref{moneqHcom24} again.\\
Let us focus on $\tilde{J}_{u}^{\mathbb{H}^1}.$ Keeping in mind the symmetry of $u,$ after a straightforward calculation, we get, because some integrals vanish by symmetry, 
 \begin{equation*}\label{moneqHcom245}\begin{split}
		&\tilde{J}_{u}^{\mathbb{H}^1}(r)=\frac{1}{r^2}\int_{D(0,r)}|\nabla_{\mathbb{G}} u(\xi)|^2\hspace{0.075cm}\Gamma (0,\xi)\hspace{0.075cm}d\xi\\
		&=\frac{c_1^2+c_2^2}{r^2}\int_{D(0,r)}((1 -6x^2-6y^2)^2+9(t^2+4x^2y^2))\Gamma (0,\xi)\hspace{0.075cm}d\xi.
		\end{split}
	\end{equation*}
	This implies
	 \begin{equation*}\label{moneqHcom2456}\begin{split}
		&\tilde{J}_{u}^{\mathbb{H}^1}(r)=(c_1^2+c_2^2)\int_{D(0,1)}((1 -6r^2(x^2+y^2)^2+9r^4(t^2+4x^2y^2))\Gamma (0,\xi)\hspace{0.075cm}d\xi.
		\end{split}
	\end{equation*}
	Therefore, $\tilde{J}_{u}^{\mathbb{H}^1}$ is monotone decreasing in a right neighborhood of $0,$ unless we are in the trivial case of $c_1=c_2=0.$ In particular, $J_{u}^{\mathbb{H}^1}$ is monotone decreasing in a right neighborhood of $0$. The same argument, which we omit for shortness, holds using the surface average integrals.
	\end{proof}
	\begin{rem}\label{lastremark}
	In principle, we may select functions, in a Carnot group $\mathbb{G},$ for which the correspondent $J_{u}^{\mathbb{G}}$ exhibits a monotone increasing behavior. The simplest case, in $\mathbb{H}^{1},$ is the following one. Let
	$v(x,y,t)=t$. Then, $v$ satisfies $\Delta_{\mathbb{H}^1}v=0$ in $\mathbb{H}^1$ by a straightforward computation. Now, since $\nabla_{\mathbb{H}^1}v=(2y,-2x)$, denoting by $\xi=(x,y,t)$, it holds, for $x_0=0,$ that
	$\left(\frac{|\nabla_{\mathbb{H}^{1}}v|}{|\nabla_{\mathbb{H}^1}|\xi||}\right)^2=4|\xi|_{\mathbb{H}^1}^2$ is sufficiently regular. Hence, we can directly compute
	$$
	\Delta_{\mathbb{H}^1}\left(\left(\frac{|\nabla_{\mathbb{H}^{1}}v|}{|\nabla_{\mathbb{H}^1}|\xi||_{\mathbb{H}^1}}\right)^2\right)= 4\hspace{0.05cm}\Delta_{\mathbb{H}^1}\left(|\xi|^2\right)=32\hspace{0.05cm}\frac{x^2+y^2}{|\xi|_{\mathbb{H}^1}^2}\geq 0.
	$$ 
	Consequently, in view of Proposition \ref{Bonfi}, Theorem \ref{maintheorem} applies, and $J_{t}^{\mathbb{H}^1}$ is monotone increasing.
	\end{rem}

	We conclude observing that, in general, the computation performed in Remark \ref{lastremark} might not be possible when  
	$$\xi\to \left(\frac{|\nabla_{\mathbb{G}} u^{\pm}(x_0\circ \xi)|}{|\nabla_{\mathbb{G}}|\xi||}\right)^2$$
	is not sufficiently regular in a neighborhood of $0$.

	\bibliographystyle{abbrv}
	\bibliography{biblio}

\def\cprime{$'$}
\begin{thebibliography}{10}

\bibitem{ACF}
H.~W. Alt, L.~A. Caffarelli, and A.~Friedman.
\newblock Variational problems with two phases and their free boundaries.
\newblock {\em Trans. Amer. Math. Soc.},
  282\hspace{0.075cm}(2):\hspace{0.075cm}431--461, 1984.

\bibitem{BLU}
A.~Bonfiglioli, E.~Lanconelli, and F.~Uguzzoni.
\newblock {\em Stratified {L}ie groups and potential theory for their
  sub-{L}aplacians}.
\newblock Springer Monographs in Mathematics. Springer, Berlin, 2007.

\bibitem{bourbaki}
N.~Bourbaki.
\newblock {\em \'{E}l\'{e}ments de math\'{e}matique. {F}asc. {XXVI}. {G}roupes
  et alg\`ebres de {L}ie. {C}hapitre {I}: {A}lg\`ebres de {L}ie}.
\newblock Actualit\'{e}s Scientifiques et Industrielles [Current Scientific and
  Industrial Topics], No. 1285. Hermann, Paris, 1971.
\newblock Seconde \'{e}dition.

\bibitem{CS}
L.~Caffarelli and S.~Salsa.
\newblock {\em A geometric approach to free boundary problems}, volume~68 of
  {\em Graduate Studies in Mathematics}.
\newblock American Mathematical Society, Providence, RI, 2005.

\bibitem{capdangar}
L.~Capogna, D.~Danielli, and N.~Garofalo.
\newblock The geometric {S}obolev embedding for vector fields and the
  isoperimetric inequality.
\newblock {\em Comm. Anal. Geom.},
  2\hspace{0.075cm}(2):\hspace{0.075cm}203--215, 1994.

\bibitem{Citti_Garofalo_Lanconelli}
G.~Citti, N.~Garofalo, and E.~Lanconelli.
\newblock Harnack's inequality for sum of squares of vector fields plus a
  potential.
\newblock {\em Amer. J. Math.},
  115\hspace{0.075cm}(3):\hspace{0.075cm}699--734, 1993.

\bibitem{Fabes_Garofalo}
E.~B. Fabes and N.~Garofalo.
\newblock Mean value properties of solutions to parabolic equations with
  variable coefficients.
\newblock {\em J. Math. Anal. Appl.},
  121\hspace{0.075cm}(2):\hspace{0.075cm}305--316, 1987.

\bibitem{FeFo2}
F.~Ferrari and N.~Forcillo.
\newblock {A counterexample to the monotone increasing behavior of an Alt
  Caffarelli-Friedman formula in the Heisenberg group}.
\newblock {\em In press, Rendiconti Lincei Matematica e Applicazioni}, arXiv:
  2203.06232v1.

\bibitem{FeFo}
F.~Ferrari and N.~Forcillo.
\newblock A new glance to the {A}lt-{C}affarelli-{F}riedman monotonicity
  formula.
\newblock {\em Math. Eng.}, 2\hspace{0.05cm}(4):\hspace{0.075cm}657--679, 2020.

\bibitem{FeLeSa}
F.~Ferrari, C.~Lederman, and S.~Salsa.
\newblock Recent results on nonlinear elliptic free boundary problems.
\newblock {\em Vietnam J. Math.},
  50\hspace{0.075cm}(4):\hspace{0.075cm}977--996, 2022.

\bibitem{Folland_Stein}
G.~B. Folland and E.~M. Stein.
\newblock {\em Hardy spaces on homogeneous groups}, volume~28 of {\em
  Mathematical Notes}.
\newblock Princeton University Press, Princeton, N.J.; University of Tokyo
  Press, Tokyo, 1982.

\bibitem{FSSC_houston}
B.~Franchi, R.~Serapioni, and F.~Serra~Cassano.
\newblock Meyers-{S}errin type theorems and relaxation of variational integrals
  depending on vector fields.
\newblock {\em Houston J. Math.},
  22\hspace{0.075cm}(4):\hspace{0.075cm}859--890, 1996.

\bibitem{FSSC_step2}
B.~Franchi, R.~Serapioni, and F.~Serra~Cassano.
\newblock Meyers-{S}errin type theorems and relaxation of variational integrals
  depending on vector fields.
\newblock {\em Houston J. Math.},
  22\hspace{0.075cm}(4):\hspace{0.075cm}859--890, 1996.

\bibitem{FSSC_CAG}
B.~Franchi, R.~Serapioni, and F.~Serra~Cassano.
\newblock On the structure of finite perimeter sets in step 2 {C}arnot groups.
\newblock {\em J. Geom. Anal.}, 13\hspace{0.075cm}(3):\hspace{0.075cm}421--466,
  2003.

\bibitem{Fulks}
W.~Fulks.
\newblock A mean value theorem for the heat equation.
\newblock {\em Proc. Amer. Math. Soc.},
  17\hspace{0.05cm}:\hspace{0.075cm}6--11, 1966.

\bibitem{GN}
N.~Garofalo and D.-M. Nhieu.
\newblock Isoperimetric and {S}obolev inequalities for
  {C}arnot-{C}arath\'eodory spaces and the existence of minimal surfaces.
\newblock {\em Comm. Pure Appl. Math.},
  49\hspace{0.075cm}(10):\hspace{0.075cm}1081--1144, 1996.

\bibitem{GT}
D.~Gilbarg and N.~S. Trudinger.
\newblock {\em Elliptic partial differential equations of second order}.
\newblock Classics in Mathematics. Springer-Verlag, Berlin, 2001.
\newblock Reprint of the 1998 edition.

\bibitem{Hormander}
L.~H\"{o}rmander.
\newblock Hypoelliptic second order differential equations.
\newblock {\em Acta Math.}, 119:\hspace{0.075cm}147--171, 1967.

\bibitem{Kupcov}
L.~P. Kupcov.
\newblock A mean value theorem and a maximum principle for a {K}olmogorov
  equation.
\newblock {\em Mat. Zametki}, 15:\hspace{0.075cm}479--489, 1974.

\bibitem{Kupcov2}
L.~P. Kupcov.
\newblock The mean value property and the maximum principle for second order
  parabolic equations.
\newblock {\em Dokl. Akad. Nauk SSSR},
  242\hspace{0.075cm}(3):\hspace{0.075cm}529--532, 1978.

\bibitem{Kupcov3}
L.~P. Kupcov.
\newblock On parabolic means.
\newblock {\em Dokl. Akad. Nauk SSSR},
  252\hspace{0.075cm}(2):\hspace{0.075cm}296--301, 1980.

\bibitem{Kuptsov}
L.~P. Kuptsov.
\newblock Property of the mean for the generalized equation of {A}. {N}.
  {K}olmogorov. {I}.
\newblock {\em Differentsial\cprime nye Uravneniya},
  19\hspace{0.075cm}(2):\hspace{0.075cm}295--304, 366--367, 1983.

\bibitem{magnani}
V.~Magnani.
\newblock Differentiability and area formula on stratified {L}ie groups.
\newblock {\em Houston J. Math.},
  27\hspace{0.05cm}(2):\hspace{0.075cm}297--323, 2001.

\bibitem{Vesely}
I.~Netuka and J.~Vesel\'{y}.
\newblock Mean value property and harmonic functions.
\newblock In {\em Classical and modern potential theory and applications
  ({C}hateau de {B}onas, 1993)}, volume 430 of {\em NATO Adv. Sci. Inst. Ser.
  C: Math. Phys. Sci.}, pages 359--398. Kluwer Acad. Publ., Dordrecht, 1994.

\bibitem{Pini_2}
B.~Pini.
\newblock Maggioranti e minoranti delle soluzioni delle equazioni paraboliche.
\newblock {\em Ann. Mat. Pura Appl. (4)}, 37:\hspace{0.075cm}249--264, 1954.

\bibitem{Pini_1}
B.~Pini.
\newblock Sulla soluzione generalizzata di {W}iener per il primo problema di
  valori al contorno nel caso parabolico.
\newblock {\em Rend. Sem. Mat. Univ. Padova}, 23:\hspace{0.075cm}422--434,
  1954.

\bibitem{stein}
E.~M. Stein.
\newblock {\em Harmonic analysis: real-variable methods, orthogonality, and
  oscillatory integrals}, volume~43 of {\em Princeton Mathematical Series}.
\newblock Princeton University Press, Princeton, NJ, 1993.
\newblock With the assistance of Timothy S. Murphy, Monographs in Harmonic
  Analysis, III.

\bibitem{Weber}
M.~Weber.
\newblock The fundamental solution of a degenerate partial differential
  equation of parabolic type.
\newblock {\em Trans. Amer. Math. Soc.}, 71:\hspace{0.075cm}24--37, 1951.

\end{thebibliography}
\end{document}